\documentclass[10pt]{article}
\usepackage[margin=1in]{geometry}
\usepackage{amsmath,amssymb,mathtools}
\usepackage{booktabs}
\usepackage{array}
\usepackage{enumitem}
\usepackage{xcolor}
\usepackage[colorlinks=true,allcolors=blue]{hyperref}
\usepackage[T1]{fontenc}
\usepackage{lmodern}
\usepackage{microtype}
\usepackage{longtable}
\usepackage{url}
\usepackage{verbatim}

\setlist{itemsep=2pt,topsep=4pt}

\newcommand{\code}[1]{\texttt{#1}}

\newcommand{\statusA}{\textbf{A}}
\newcommand{\statusB}{\textbf{B}}
\newcommand{\statusC}{\textbf{C}}
\newcommand{\statusD}{\textbf{D}}

\title{A Lean 4 Verification Report for\\[4pt]
\emph{Completing the Arakawa--Moreau Conjecture on Maximal Ideals of Affine Vertex Algebras}}
\author{Sihai Jin\\
Department of Mathematics, Sichuan University\\
\texttt{jinsihai@stu.scu.edu.cn}}
\date{}

\begin{document}
\maketitle

\begin{abstract}
We report a Lean 4 formal audit accompanying the paper
\emph{Completing the Arakawa--Moreau Conjecture on Maximal Ideals of Affine Vertex Algebras}
(arXiv:2607.25249).  In a single Lean source, the development kernel-checks the paper-local deductive architecture used for the new maximal-ideal cases: affine-kernel detection and lifting, extremal PBW-line deductions, finite Ramond--Casimir arithmetic and case enumerations, reduced-simplicity contradiction schemes, and the case-level affine conclusions for the level $-1$ $D_\ell$ family and the positive-parameter cases of $D_4,E_6,E_7,E_8$.  It also audits the paper's alternative $D_\ell$ level $-2$ argument and the separate rank-reduction application of Section~8.  Previously published exceptional $n=0$ cases are not claimed as newly formalized results of this bundle.  The released project builds successfully with the pinned Lean/Mathlib environment and contains live \code{\#print axioms} queries for the principal final wrappers.

The development does not reconstruct affine vertex algebras, minimal $W$-algebras, BRST/DS reduction, Ramond Zhu theory, or Li spectral flow from foundational definitions in Mathlib.  Those ingredients, together with the concrete realization of case-specific representation-theoretic objects, are exposed as theorem parameters and semantic interfaces.  Accordingly, the precise verification claim is a kernel-checked deduction of the paper-local proof from explicitly stated VOA/BRST/DS boundary inputs, rather than a from-scratch formalization of the ambient representation theory.
\end{abstract}

\section{Purpose and relation to the original paper}

The original paper proves the remaining cases of Arakawa--Moreau Conjecture~1 concerning explicit maximal ideals of negative-level universal affine vertex algebras of types $D$ and $E$.  Its main theorem states that the prescribed singular vectors generate the maximal graded ideals in every negative level covered by the conjecture.  In the $D_\ell$ series, the new result is the level $-1$ statement for every $\ell\ge 5$, while the level $-2$ statement is also reproved by a different reduction-theoretic method.  For the exceptional types, the new results are the cases with positive parameter $n$ for $D_4,E_6,E_7,E_8$; together with previously known $n=0$ cases, these complete the conjecture.

The released source is available both as arXiv ancillary material and in the public GitHub repository
\url{https://github.com/jshemail12345-debug/ArakawaMoreauFormalAudit}; the version corresponding to this report is frozen as release
\url{https://github.com/jshemail12345-debug/ArakawaMoreauFormalAudit/releases/tag/v1.0.0}.

The accompanying Lean source is a single-file formal-audit bundle.  Its purpose is not to recreate the full background theory from basic definitions.  Instead, the source separates the argument into two layers:
\begin{enumerate}[label=(\roman*)]
  \item paper-local arithmetic, finite combinatorics, abstract PBW saturation, reduced-simplicity logic, and affine lifting checked directly by the Lean kernel; and
  \item foundational or literature-level VOA/BRST/DS facts represented by explicit hypotheses or structured interfaces.
\end{enumerate}

The principal released source file is
\begin{center}
\code{ArakawaMoreauFormalAudit.lean}.
\end{center}
The project is designed so that a successful \code{lake build} checks all declarations and executes the live axiom-dependency audit at the end of the source.

\subsection*{Verification statement}
\noindent\fbox{\begin{minipage}{0.94\textwidth}
\textbf{Claim certified by the released bundle.}
With the explicitly listed foundational and case-realization hypotheses in the Lean theorem signatures, the paper-local deductions implemented in \code{ArakawaMoreauFormalAudit.lean} are accepted by the Lean kernel.  The frozen GitHub release \code{v1.0.0} completed \code{lake build} successfully before publication.  No project-local declaration of the form \code{axiom ...} and no proof placeholder introduced by \code{sorry} or \code{admit} occurs in the frozen source.
\end{minipage}}

This statement is deliberately narrower than ``the whole paper has been formalized from definitions.''  The distinction is part of the formal audit and is repeated in Section~\ref{sec:claim-boundary}.

\section{Scope of the formalization}

\subsection{Kernel-checked internal components}

The source records the following principal internally checked components.
\begin{enumerate}
  \item The abstract affine-kernel contradiction underlying Proposition~2.10 and Theorem~2.11: exactness forces vanishing of the reduced kernel, while nonvanishing/detection rules out a nonzero affine kernel.
  \item Lowest nonzero kernel-degree extraction by well-ordering and the negative-level affine-wall inequality used in the Lean version of Lemma~2.9.
  \item The generic extremal PBW saturation argument corresponding to the logical core of Proposition~2.4.
  \item Conversion of a nonzero reduced singular-vector class lying on an extremal PBW line into a nonzero scalar multiple of the corresponding current power.
  \item Invariance of generated ideals under nonzero rescaling and the exact quotient-algebra deduction corresponding to the logic of Proposition~2.3 once the external exactness input is supplied.
  \item Weighted-average and Casimir-gap lemmas used to turn pointwise constituent bounds into the contradiction required for reduced simplicity.
  \item Ground-energy forcing and extremal-state-to-vacuum logic in the abstract Ramond interface.
  \item The finite arithmetic, Casimir bounds, and finite enumerations encoded for the $D_\ell$, $D_4$, $E_6$, $E_7$, and $E_8$ cases.
  \item Case-level maximality wrappers corresponding to the paper's final affine statements in Sections~3--7.
  \item The Section~8 sandwich/rank-reduction argument and its graded affine-lifting step.
  \item Scalar-generic exact-image wrappers for the $D_\ell$, $D_4$, $E_6$, $E_7$, $E_8$, and rank-reduction exact-image statements.
\end{enumerate}

\subsection{Coverage at a glance}

\begin{center}
\small
\begin{tabular}{>{\raggedright\arraybackslash}p{0.18\textwidth} >{\raggedright\arraybackslash}p{0.34\textwidth} >{\raggedright\arraybackslash}p{0.38\textwidth}}
\toprule
Paper part & Kernel-checked layer & Main boundary layer\\
\midrule
Section 2 & Kernel/detection logic, PBW saturation, weighted-gap deductions & DS exactness/nonvanishing, PBW realization, Ramond/Li semantics\\
Section 3 ($D_\ell$) & Uniform wall, exponent, Casimir and gap arithmetic; final lifting & Concrete reduced classes, lowest Ramond constituents\\
Section 4 ($D_4$) & Triality PBW defect, finite $A_1^3$ enumeration, trace/gap logic & Concrete BRST/Ramond realization\\
Sections 5--7 & Type-specific finite Casimir tables, pointwise-to-average gaps, final lifting & Concrete $E_6,E_7,E_8$ constituent extraction and DS realization\\
Section 8 & Scalar-line extraction, sandwich/bijectivity and induction bookkeeping & Concrete rank-reduction BRST/generation facts\\
\bottomrule
\end{tabular}
\end{center}

\subsection{Foundational and literature boundary}

The formalization deliberately leaves the following material as explicit inputs rather than rebuilding it from definitions.
\begin{enumerate}
  \item Construction and basic algebraic properties of the relevant universal and simple affine vertex algebras and universal/simple minimal $W$-algebras.
  \item Exactness of minimal DS reduction on the relevant affine category.
  \item The irreducible highest-weight DS nonvanishing theorem away from the affine wall.
  \item Reduction of affine Verma modules to $W$-Verma modules and the cyclicity consequences for reduced singular submodules.
  \item The standard PBW spanning theorem for universal minimal $W$-algebras and the identification of concrete current/$G$/Virasoro generators with the abstract PBW atoms used by the Lean interface.
  \item The universal Ramond Zhu relation, Li spectral-flow formulas, and the representation-theoretic extraction of finite-dimensional Zhu quotients and constituents from a lowest Ramond space.
  \item Concrete bridges identifying the abstract Lean types and predicates with the actual type-dependent $W$-algebra states, modules, ideals, Casimir constituents, and highest-weight data in the paper.
  \item Previously published cases of Arakawa--Moreau Conjecture~1 when they are used as literature rather than reproved inside the paper's alternative argument.
\end{enumerate}

These are theorem parameters rather than custom Lean axioms.  Consequently, \code{\#print axioms} certifies the dependency of the formal deduction on Lean's logical foundations, but it does not by itself certify the mathematical truth of the external theorem parameters.

\section{Detailed audit convention}

For the section-by-section discussion below, it is useful to distinguish four levels of formalization.  The classification is descriptive; it does not change the Lean trust model.
\begin{description}[leftmargin=1.8cm,style=nextline]
  \item[\statusA: internal Lean deduction.] Concrete arithmetic, finite combinatorics, linear or order-theoretic deductions, abstract PBW saturation, ideal/kernel logic, and other conclusions are proved inside Lean from lower-level hypotheses.
  \item[\statusB: paper-formula boundary.] A formula or numerical identity derived in the mathematical paper from representation theory is supplied as an input, while a substantial downstream simplification, inequality, finite elimination, or contradiction is checked in Lean.
  \item[\statusC: theorem boundary.] A foundational or literature-level theorem---for example DS exactness/nonvanishing, cyclicity of a reduced singular submodule, the universal Ramond--Zhu relation, or a simplicity/nonvanishing fact for a standard reduction---is supplied semantically rather than reconstructed from definitions.
  \item[\statusD: realization bridge.] A hypothesis connects the actual paper object (a VOA state, BRST class, PBW generator, Ramond constituent, ideal, or quotient) with the abstract type or predicate on which the kernel-checked calculation is performed.
\end{description}
A result may involve more than one level.  For example, ``A+C/D'' means that the final deduction is internal once theorem-level and realization-level interfaces are provided.

\section{Detailed audit of the common framework (Section 2)}

Section~2 of the paper isolates the reusable mechanism behind all type-dependent arguments.  The Lean bundle mirrors this separation and contains the most reusable part of the formal audit.

\subsection{Affine detection and lifting}

The core Lean theorem \path{kernel_trivial_of_exactness_detection} formalizes the following abstract implication.  If exact reduction implies that the reduced kernel is zero, and if every nonzero affine kernel produces a nonzero reduced kernel, then the original quotient map is injective.  The declarations \path{affine_lifting_injective}, \path{affine_lifting_bijective}, \path{candidate_eq_radical_of_injective_quotient_map}, and \path{minimal_reduction_maximality_core} package the remaining kernel and ideal-equality logic.

At the paper-specific interface level, \path{lemma_2_9_detection_from_graded_kernel_literature}, \path{proposition_2_10_from_graded_kernel_literature}, and \path{theorem_2_11_from_graded_kernel_literature} make the trust boundary explicit.  Lean internalizes the lowest-nonzero-degree contradiction and the formal transfer from injectivity to equality with the radical.  The existence of the appropriate affine highest-weight witness, DS exactness, and the irreducible DS nonvanishing theorem remain Level-C/D inputs.

\subsection{Extremal PBW line and exact-image logic}

The paper's Proposition~2.4 is represented by an abstract PBW saturation interface.  Lean checks the combinatorial implication that equality in the charge-versus-degree estimate forces every atom of a saturated monomial to be the unique extremal current atom.  It also checks the passage from one-dimensionality of the extremal line and nonvanishing of a reduced class to a nonzero scalar multiple of the current power.  The scalar-generic layer is formulated over an arbitrary field, rather than hard-coding complex scalars.

The realization of the actual universal minimal $W$-algebra PBW basis by the abstract atoms, and the statement that the concrete BRST class lies in the required extremal space, are Level-D interfaces.  The DS nonvanishing/cyclicity statements are Level C.  Thus the final scalar-line conclusion is kernel-checked, but the BRST complex itself is not reconstructed.

\subsection{Ramond--Zhu finiteness and Casimir-gap mechanism}

The source contains generic finite-cyclic-submodule and simple-quotient constructions together with weighted-average bounds.  Theorems such as \path{weighted_average_le_max}, \path{casimir_gap_forces_ground_energy}, \path{ground_energy_vacuum_contradiction}, and \path{casimir_gap_simplicity_core} isolate the logical content of Proposition~2.7.  Later refinements prove the same conclusion from pointwise constituent bounds, rather than assuming a pre-averaged estimate.

The representation-theoretic extraction of a finite Ramond Zhu module from a lowest ideal energy, the universal Ramond relation, and the case-specific identification of the extremal constituent with the vacuum state remain explicit interfaces.  Given those data, Lean checks the weighted-average contradiction and the ideal-equals-whole-algebra conclusion.

\subsection{Section-2 status table}

{\footnotesize
\begin{longtable}{>{\raggedright\arraybackslash}p{0.17\textwidth} >{\raggedright\arraybackslash}p{0.12\textwidth} >{\raggedright\arraybackslash}p{0.27\textwidth} >{\raggedright\arraybackslash}p{0.34\textwidth}}
\toprule
Paper step & Level & Representative Lean content & Verification boundary\\
\midrule
\endfirsthead
\toprule
Paper step & Level & Representative Lean content & Verification boundary\\
\midrule
\endhead
Lemma 2.1 & A+C/D & graded-kernel/lowest-degree infrastructure & Finite graded pieces and membership in the DS exactness category are not rebuilt as actual affine modules.\\
Lemma 2.2 & C/D & consumed by exact-image nonvanishing interfaces & Filtered BRST spectral-sequence and Slodowy-slice restriction theorem remain external.\\
Proposition 2.3 & A+C & quotient/kernel logic and generated-ideal transport & Exactness and cyclicity of the reduced singular submodules are supplied.\\
Proposition 2.4 & A+D & extremal PBW saturation and scalar-line deduction & Concrete $W$-algebra PBW generators and class-membership bridges are supplied.\\
Proposition 2.5 & A+B/C & weighted trace arithmetic and case specializations & Universal Ramond--Zhu relation and contraction origin are not foundationally formalized.\\
Lemma 2.6 & A+C/D & finite cyclic submodule/simple quotient infrastructure & Semantic interpretation as the lowest Ramond Zhu space is an interface.\\
Proposition 2.7 & A+C/D & \path{casimir_gap_simplicity_core} and pointwise refinements & Constituents, trace identity, and vacuum-state realization are supplied casewise.\\
Lemma 2.9 & A+C/D & \path{lemma_2_9_detection_from_graded_kernel_literature} & Affine highest-weight/DS nonvanishing theorem is external.\\
Proposition 2.10 & A+C & \path{proposition_2_10_from_graded_kernel_literature} & DS exactness/nonvanishing bridge supplied.\\
Theorem 2.11 & A+C & \path{theorem_2_11_from_graded_kernel_literature} & Final kernel-to-radical equality is internal from the displayed interfaces.\\
Theorem 2.13 & A+C/D & \path{ds_ramond_casimir_maximality_core} and case wrappers & Packages the reduced-simplicity and affine-lifting architecture without reconstructing VOA foundations.\\
\bottomrule
\end{longtable}
}

\section{Detailed audit of the $D_\ell$ series (Section 3)}

The $D_\ell$ source is the most uniform infinite-rank part of the development.  The Lean file checks the rank arithmetic, wall inequalities, Casimir bounds, and gap formulas uniformly in $\ell\ge5$, and then reuses the common Section~2 logic.

\subsection{Level $-1$: finite arithmetic and wall data}

The definitions \code{mOfRank} and \code{KOfRank} encode $m=\ell-2$ and $K=m-1=\ell-3$.  Lean proves the two generator exponents and the affine-wall inequalities used by the paper's nonvanishing argument.  It then encodes the allowed $A_1$ and surviving $D_m$ constituent types and proves the paper's pointwise Casimir bounds.  The declarations \code{weightedCasimir\_le}, \code{first\_gap\_formula}, \code{first\_casimir\_gap}, and \code{positive\_energy\_casimir\_gap} verify the numerical heart of Theorem~3.21 uniformly in the rank.

The mathematical input that the paper's nilpotence and weight-coset arguments reduce the actual Ramond constituents to the encoded allowed types is a Level-C/D bridge.  Lean does not silently assume the final average inequality: it proves the weighted bound from the encoded pointwise bounds.

\subsection{Exact reduced singular vectors and reduced quotient}

For Proposition~3.6 the wrapper \path{DSeries.proposition_3_6_exact_images_over_field} checks, over an arbitrary coefficient field, that each nonzero BRST class lying in the saturated extremal PBW space is a nonzero scalar multiple of the appropriate current power.  The nonvanishing, cyclicity, PBW-spanning, and concrete membership statements appear explicitly in the signature.  This matches the paper's logical use of Lemmas~3.3--3.5 without claiming a foundational BRST implementation.

The exact reduction of the candidate quotient (Proposition~3.8) is represented by the common quotient logic from Proposition~2.3 together with the case-specific exact-image interface.  Associated-variety and lisse inputs remain external where used to obtain a finite Ramond Zhu algebra.

\subsection{Reduced simplicity and affine maximality at level $-1$}

Theorem~3.21 is represented in several layers.  \path{DSeries.theorem_3_21_simplicity_core} exposes the pure contradiction pattern; \path{DSeries.theorem_3_21_reduced_simplicity_certificate} packages the reduced-simplicity certificate.  The subsequent wrappers \path{DSeries.theorem_3_23_gap_to_affine} and \path{DSeries.theorem_3_23_gap_to_affine_graded_literature} connect that certificate to the graded affine-lifting interface.  Thus the strict Casimir gap and kernel contradiction are Level A, whereas the realization of the lowest Ramond space and the DS exactness/detection theorem remain Level C/D.

\subsection{Endpoint level $-2$}

The endpoint is treated separately, just as in the paper.  Lean defines the endpoint trace function and Casimir maximum and proves the exact first-shift gap.  The scalar-generic wrapper \path{DSeries.proposition_3_27_exact_images_over_field} handles Proposition~3.27.  The declarations \path{DSeries.theorem_3_33_simplicity_core} and \path{DSeries.theorem_3_34_gap_to_affine_graded_literature} formalize the endpoint reduced-simplicity and affine-lifting deductions.  The collapsing identification of the reduced algebra with the surviving $A_1$ factor is a representation-theoretic input rather than a foundational construction.

\subsection{$D_\ell$ status table}

{\footnotesize
\begin{longtable}{>{\raggedright\arraybackslash}p{0.17\textwidth} >{\raggedright\arraybackslash}p{0.12\textwidth} >{\raggedright\arraybackslash}p{0.28\textwidth} >{\raggedright\arraybackslash}p{0.33\textwidth}}
\toprule
Paper result & Level & Lean declaration / package & Status\\
\midrule
\endfirsthead
\toprule
Paper result & Level & Lean declaration / package & Status\\
\midrule
\endhead
Lemmas 3.3--3.5 / Prop. 3.6 & A+C/D & \path{DSeries.proposition_3_6_exact_images_over_field} & Wall arithmetic and extremal-line deduction internal; DS nonvanishing/cyclicity and concrete PBW realization external.\\
Proposition 3.8 & A+C & common exact-quotient package & Quotient logic internal; exactness, associated variety, and lisse facts supplied.\\
Lemmas 3.10--3.16 & A+B/C/D & \path{a1Casimir_le}, \path{dCasimir_le}, \path{weightedCasimir_le} & Numerical bounds internal after the paper's constituent restriction is encoded.\\
Proposition 3.12 & A+B/C & trace/gap arithmetic package & Universal trace relation/contraction origin supplied; numerical consequences checked.\\
Props. 3.17--3.20 & A+C/D & ground-energy interfaces plus gap package & Spectral lattice/bottom-space realization is not rebuilt; consequences used by the gap argument are explicit.\\
Theorem 3.21 & A+C/D & \path{DSeries.theorem_3_21_reduced_simplicity_certificate} & Reduced-simplicity contradiction kernel-checked.\\
Theorem 3.23 & A+C & \path{DSeries.theorem_3_23_gap_to_affine_graded_literature} & Final maximal-ideal equality checked from graded DS interfaces.\\
Proposition 3.27 & A+C/D & \path{DSeries.proposition_3_27_exact_images_over_field} & Endpoint scalar-line conclusion internal from explicit nonvanishing/PBW inputs.\\
Props. 3.31--3.32 & A+B/C/D & endpoint trace and ground package & Endpoint numerical trace gap internal; spectral interpretation external.\\
Theorem 3.33 & A+C/D & \path{DSeries.theorem_3_33_reduced_simplicity_certificate} & Endpoint reduced simplicity checked from stated inputs.\\
Theorem 3.34 & A+C & \path{DSeries.theorem_3_34_gap_to_affine_graded_literature} & Endpoint affine lifting/kernel equality internal.\\
\bottomrule
\end{longtable}
}

\section{Detailed audit of the triality case $D_4$ (Section 4)}

The $D_4$ case has the densest case-specific internal combinatorics in the file.  The source explicitly models the three triality factors, their allowed $A_1^3$ data, the PBW charge/degree defect, the Ramond trace average, and the lowest-energy contradiction.

\subsection{Exact-image and extremal-line calculation}

Lean proves the uniform triality wall value and its negativity, develops an atom-by-atom charge-versus-doubled-degree defect, and proves that zero total defect forces every atom to be the distinguished extremal current.  The declarations \path{extremal_monomial_classification}, \path{extremal_monomial_eq_pair}, \path{extremal_space_lies_in_square_line}, and \path{extremal_space_one_dimensional_line} make the extremal-square line explicit in the abstract PBW model.  The release-level statement is \path{D4.proposition_4_3_exact_images_over_field}.

The actual BRST class, its nonvanishing, and its placement in the concrete extremal $W$-space are Level-C/D inputs at the final generic wrapper.  Earlier local arithmetic such as the wall evaluation and reduced degree/weight profile is checked internally.

\subsection{Finite $A_1^3$ types and the Casimir trace}

The structure \code{AllowedTriple} encodes the finite set of allowed highest-weight triples.  Lean defines each $A_1$ Casimir, the total Casimir, the maximal triple, and $C_{\max}$, and proves the uniform pointwise bound.  The source also checks the specialization of the traced Ramond relation and its identification with the affine function $A(h)$ through \path{ramond_trace_identity_from_uniform_formula}, \path{ramond_trace_identity}, and \path{trace_identity_matches_A}.

The universal Ramond relation and the semantic fact that the actual lowest ideal space yields one of the encoded finite constituents remain boundary inputs.  Once they are supplied, the finite arithmetic is kernel-checked.

\subsection{Spectral bottom, reduced simplicity, and affine lifting}

The PBW contribution model proves that all Ramond energies lie above the vacuum and in the required lattice.  Separate finite arguments derive the bottom degree--weight restrictions and rule out the competing ground types.  The file then builds finite cyclic Zhu submodules and simple quotients and packages the conclusion as \path{D4.d4_reduced_simplicity_certificate}.  The final paper-level wrapper \path{D4.d4_theorem_4_14_graded_literature} combines this with graded affine detection/lifting.

\subsection{$D_4$ status table}

{\footnotesize
\begin{longtable}{>{\raggedright\arraybackslash}p{0.17\textwidth} >{\raggedright\arraybackslash}p{0.12\textwidth} >{\raggedright\arraybackslash}p{0.28\textwidth} >{\raggedright\arraybackslash}p{0.33\textwidth}}
\toprule
Paper result & Level & Lean declaration / package & Status\\
\midrule
\endfirsthead
\toprule
Paper result & Level & Lean declaration / package & Status\\
\midrule
\endhead
Lemma 4.1 & A+C/D & \path{lemma_4_1_core} and wall/profile arithmetic & Wall and cyclic-generator logic internal; DS nonvanishing theorem supplied.\\
Lemma 4.2 / Prop. 4.3 & A+C/D & extremal PBW package; \path{D4.proposition_4_3_exact_images_over_field} & Extremal-square uniqueness and scalar-line conclusion internal from PBW realization inputs.\\
Proposition 4.4 & A+C & common exact quotient logic & Exactness/cyclicity and concrete quotient interpretation supplied.\\
Lemmas 4.5--4.7 & A+B/C/D & \path{AllowedTriple}, Casimir/contraction package & Finite triple/Casimir arithmetic internal; representation-theoretic allowed-type and universal relation inputs explicit.\\
Proposition 4.8 & A+B/C & \path{D4.proposition_4_8_from_traced_relation}, \path{D4.proposition_4_8_matches_A} & Trace specialization and simplification kernel-checked.\\
Props. 4.10--4.12 & A+C/D & PBW spectral lattice and bottom degree--weight package & Concrete Li-twist/Ramond interpretation is an interface.\\
Theorem 4.13 & A+C/D & \path{D4.d4_reduced_simplicity_certificate} & Finite gap/exclusion and vacuum contradiction internal.\\
Theorem 4.14 & A+C & \path{D4.d4_theorem_4_14_graded_literature} & Affine kernel vanishing and candidate-ideal equality checked from graded DS inputs.\\
\bottomrule
\end{longtable}
}

\section{Detailed audit of type $E_6$ (Section 5)}

The $E_6$ formalization follows the same architecture with two allowed values $n=1,2$.  The source checks the level/exponent arithmetic, the specialized trace functions, the global and second Casimir bounds, the strict positive-shift gap, and the ground-extremal gap.

The exact-image conclusion of Proposition~5.7 is represented by \path{E6.proposition_5_7_exact_images_over_field}.  The reduced-simplicity theorem has both certificate and pointwise forms, notably \path{E6.theorem_5_18_reduced_simplicity_certificate} and \path{E6.theorem_5_18_reduced_simplicity_pointwise}.  The latter derives the averaged contradiction from pointwise constituent bounds, so the final average bound is not simply inserted as a hypothesis.  The affine conclusion is \path{E6.theorem_5_19_pointwise_graded_literature}.

{\footnotesize
\begin{longtable}{>{\raggedright\arraybackslash}p{0.17\textwidth} >{\raggedright\arraybackslash}p{0.12\textwidth} >{\raggedright\arraybackslash}p{0.28\textwidth} >{\raggedright\arraybackslash}p{0.33\textwidth}}
\toprule
Paper result & Level & Lean declaration / package & Status\\
\midrule
\endfirsthead
\toprule
Paper result & Level & Lean declaration / package & Status\\
\midrule
\endhead
Theorem 5.6 / Prop. 5.7 & A+C/D & level/wall arithmetic; \path{E6.proposition_5_7_exact_images_over_field} & Singular-vector theorem and BRST realization supplied; exact scalar current-power deduction checked.\\
Props. 5.12--5.14 & A+B/C/D & \path{E6.trace_specialization_one}, \path{E6.trace_specialization_two}, gap package & Exact reduction/lisse and universal Ramond relation external; specialized arithmetic internal.\\
Lemma 5.15 & A+B/D & $C_{\max}$ and $C_{\mathrm{next}}$ bounds encoded pointwise & Finite weight classification bridge supplied; downstream inequalities checked.\\
Proposition 5.17 & A+C/D & ground-energy/extremal-state interface plus internal gap logic & Actual Li twist and bottom-space representation are not reconstructed.\\
Theorem 5.18 & A+C/D & \path{E6.theorem_5_18_reduced_simplicity_pointwise} & Pointwise Casimir data imply reduced simplicity in Lean.\\
Theorem 5.19 & A+C & \path{E6.theorem_5_19_pointwise_graded_literature} & Final affine lifting/maximal-ideal conclusion kernel-checked.\\
\bottomrule
\end{longtable}
}

\section{Detailed audit of type $E_7$ (Section 6)}

For $E_7$ the allowed values are $n=1,2,3$.  Lean proves the negative wall values, the specialized trace average, the first positive-shift gap, and the ground-extremal gap.  The finite $D_6$ Casimir data are encoded so that the weighted-average elimination is performed by the kernel.

The exact-image wrapper is \path{E7.proposition_6_4_exact_images_over_field}.  The reduced-simplicity pointwise theorem is \path{E7.theorem_6_18_reduced_simplicity_pointwise}, and the final affine wrapper is \path{E7.theorem_6_19_pointwise_graded_literature}.  The low-degree representation-theoretic restrictions used to exclude competing $D_6$ constituents are exposed as boundary data; Lean checks the finite implication from these pointwise restrictions to the extremal constituent and then to simplicity.

{\footnotesize
\begin{longtable}{>{\raggedright\arraybackslash}p{0.17\textwidth} >{\raggedright\arraybackslash}p{0.12\textwidth} >{\raggedright\arraybackslash}p{0.28\textwidth} >{\raggedright\arraybackslash}p{0.33\textwidth}}
\toprule
Paper result & Level & Lean declaration / package & Status\\
\midrule
\endfirsthead
\toprule
Paper result & Level & Lean declaration / package & Status\\
\midrule
\endhead
Props. 6.1--6.4 & A+C/D & wall/exponent arithmetic; \path{E7.proposition_6_4_exact_images_over_field} & Nonvanishing/cyclicity and concrete PBW realization supplied; scalar-line conclusion internal.\\
Props. 6.5--6.10 & A+B/C/D & contraction and trace-specialization package & Exact reduction/lisse and universal Ramond relation supplied; finite trace arithmetic checked.\\
Props. 6.12--6.16 & A+C/D & spectral/ground interfaces with internal degree--weight consequences & Li-twist and semantic constituent realization remain external.\\
Lemma 6.17 & A+B/D & $D_6$ finite Casimir tables and pointwise bounds & Finite case arithmetic checked; bridge to actual Zhu constituents explicit.\\
Theorem 6.18 & A+C/D & \path{E7.theorem_6_18_reduced_simplicity_pointwise} & Ground/extremal forcing and ideal contradiction checked.\\
Theorem 6.19 & A+C & \path{E7.theorem_6_19_pointwise_graded_literature} & Final affine maximality checked from graded DS interfaces.\\
\bottomrule
\end{longtable}
}

\section{Detailed audit of type $E_8$ (Section 7)}

The $E_8$ case contains the largest finite parameter range, $1\le n\le5$.  Lean checks the level values and wall negativity, the uniform trace function, the global Casimir maximum, the second-largest bound used away from the level-one exception, and the strict inequalities at every allowed $n$.  In particular, declarations such as \path{first_gap_formula}, \path{first_gap}, \path{ground_minus_second_nonlevelone}, \path{ground_extremal_gap}, and \path{positive_shift_gap} make the finite numerical gap transparent.

The exact-image statement is \path{E8.corollary_7_6_exact_images_over_field}.  The pointwise reduced-simplicity wrapper \path{E8.theorem_7_23_reduced_simplicity_pointwise} separates the numerical contradiction from the representation-theoretic extraction of actual $E_7$ constituents.  The final affine statement is \path{E8.theorem_7_24_pointwise_graded_literature}.

{\footnotesize
\begin{longtable}{>{\raggedright\arraybackslash}p{0.17\textwidth} >{\raggedright\arraybackslash}p{0.12\textwidth} >{\raggedright\arraybackslash}p{0.28\textwidth} >{\raggedright\arraybackslash}p{0.33\textwidth}}
\toprule
Paper result & Level & Lean declaration / package & Status\\
\midrule
\endfirsthead
\toprule
Paper result & Level & Lean declaration / package & Status\\
\midrule
\endhead
Theorem 7.2 / Cor. 7.6 & A+C/D & wall/exponent arithmetic; \path{E8.corollary_7_6_exact_images_over_field} & Singular-vector/BRST facts supplied; exact nonzero scalar line checked.\\
Props. 7.7--7.14 & A+B/C/D & reduction and trace-specialization package & Exactness/lisse/universal Ramond relation external; numerical specialization internal.\\
Lemmas 7.16--7.17 / Thm. 7.18 & A+B/D & encoded finite $E_7$ Casimir bounds and gap functions & Finite inequalities checked; bridge identifying actual allowed weights explicit.\\
Corollary 7.19 & A+B/D & ground second-Casimir exclusion package & Exceptional $n=1$ semantic input separated from the general finite inequality.\\
Proposition 7.22 & A+C/D & ground-energy interface and internal shift arithmetic & Concrete Li spectral flow not foundationally rebuilt.\\
Theorem 7.23 & A+C/D & \path{E8.theorem_7_23_reduced_simplicity_pointwise} & Positive-shift and ground-extremal contradictions kernel-checked.\\
Theorem 7.24 & A+C & \path{E8.theorem_7_24_pointwise_graded_literature} & Final affine lifting/maximal-ideal equality kernel-checked.\\
\bottomrule
\end{longtable}
}

\section{Detailed audit of the rank-reduction application (Section 8)}

Section~8 is logically separate from the new cases of Conjecture~1, but it is a useful test that the common lifting mechanism is reusable.  Lean checks the surviving-level arithmetic, the nonzero coefficient used to eliminate the conformal generator, and a generic sandwich lemma: if two composable maps are surjective and their composite is bijective, then the intermediate map is bijective.  The declarations \path{sandwich_bijective}, \path{induction_from_two}, and \path{theorem_8_5_affine_core} isolate these deductions.

The exact-image scalar-line component of Proposition~8.1 is \path{RankReduction.proposition_8_1_exact_images_over_field}.  The later wrapper \path{RankReduction.theorem_8_5_induction_step_graded_literature} connects the rank-reduction isomorphism to the same graded affine-lifting interface used earlier.  The concrete BRST computation identifying the reduced generators with the previous-rank generators, and the VOA generation/quadratic-relation facts used to construct the rank-reduction surjection, remain theorem/realization inputs where not reconstructed.

{\footnotesize
\begin{longtable}{>{\raggedright\arraybackslash}p{0.17\textwidth} >{\raggedright\arraybackslash}p{0.12\textwidth} >{\raggedright\arraybackslash}p{0.28\textwidth} >{\raggedright\arraybackslash}p{0.33\textwidth}}
\toprule
Paper result & Level & Lean declaration / package & Status\\
\midrule
\endfirsthead
\toprule
Paper result & Level & Lean declaration / package & Status\\
\midrule
\endhead
Proposition 8.1 & A+C/D & \path{RankReduction.proposition_8_1_exact_images_over_field} & Nonzero-on-line scalar extraction internal; concrete BRST-image line statements supplied.\\
Lemma 8.2 & A+C/D & \path{RankReduction.lemma_8_2_generation_core} & Abstract generation consequence checked; actual $W$-generator vanishing/VOA generation facts supplied.\\
Lemma 8.3 & A+B/C/D & \path{RankReduction.lemma_8_3_quadratic_relation_core} & Algebraic elimination checked from the supplied $G$--$G$/quadratic relation data.\\
Proposition 8.4 & A+C/D & \path{RankReduction.proposition_8_4_factor_core} & Factor/surjection logic internal; concrete current quotient realization external.\\
Theorem 8.5 & A+C & \path{RankReduction.theorem_8_5_induction_step_graded_literature} & Sandwich/bijectivity and affine lifting internal; base case and concrete rank-reduction bridge supplied.\\
\bottomrule
\end{longtable}
}

\section{Relation to the paper's Theorem 1.1}

The released Lean source does not introduce a single declaration named \code{theorem\_1\_1}.  Instead, it kernel-checks the new case-level maximality deductions separately: Theorem~3.23 for $D_\ell$ at level $-1$, Theorem~4.14 for $D_4$ at level $-1$, Theorem~5.19 for $E_6$, Theorem~6.19 for $E_7$, and Theorem~7.24 for $E_8$, together with the alternative endpoint Theorem~3.34 at level $-2$.  These are exactly the new case statements, together with the paper's alternative proof of the already known $D_\ell$ level $-2$ statement, whose mathematical combination with the previously published exceptional $n=0$ cases yields Theorem~1.1 of the paper.  The Lean bundle does not separately reconstruct those previously published exceptional boundary cases.  Section~8 is an additional application and is not part of Theorem~1.1.

This organization is intentional.  It keeps the formal trust boundary visible at the level where each type-dependent Ramond/DS input enters, rather than hiding all case-specific interfaces behind a single global proposition.

\section{Paper-to-Lean verification map}

Before freezing this version, the numbered statements in the table below were cross-checked against the statement labels in the paper, and every explicitly named Lean declaration was cross-checked against the frozen source.  This is a source-to-report consistency check; it does not enlarge the mathematical trust boundary.  Generic descriptions such as ``spectral-bottom package'' are used only where the source deliberately spreads a paper argument across several auxiliary declarations rather than exposing a single declaration with the paper's number.

Table~\ref{tab:paper-lean-map} gives a more granular map from the numbered statements used in the paper to the corresponding Lean declarations or proof packages.  The map is deliberately conservative: when the Lean declaration verifies only the paper-local deduction after a representation-theoretic input has been supplied, the final column says so explicitly.
For readability, declarations are written with their case namespace (for example, \code{E7.} or \code{DSeries.}); all such names live under the outer namespace \code{ArakawaMoreau}.

{\scriptsize
\begin{longtable}{>{\raggedright\arraybackslash}p{0.16\textwidth} >{\raggedright\arraybackslash}p{0.43\textwidth} >{\raggedright\arraybackslash}p{0.31\textwidth}}
\caption{Expanded paper-to-Lean verification map.}\label{tab:paper-lean-map}\\
\toprule
Paper result & Lean declaration / package & Verification boundary\\
\midrule
\endfirsthead
\toprule
Paper result & Lean declaration / package & Verification boundary\\
\midrule
\endhead
Proposition 2.3 & \path{exact_reduction_family_quotient_of_scaling}; \path{exact_reduction_of_case_family} & Exact-sequence/quotient and rescaling logic internal; DS exactness and identification of the reduced ideal range are inputs.\\
Proposition 2.4 & \path{ExtremalPBW.*}; \path{ScalarGeneric.saturated_space_lies_in_currentPower_line} & Charge-versus-degree saturation and one-line consequence internal; concrete universal-$W$ PBW realization external.\\
Proposition 2.5 & weighted-average / traced-Casimir core used by the case packages & Trace manipulation and rational simplification internal once the universal Ramond relation and contraction coefficients are supplied.\\
Lemma 2.6 / Proposition 2.7 & \path{reduced_simplicity_from_pointwise_Ramond_data}; case certificates & Lowest-shift contradiction and extremal-state-to-vacuum forcing internal; existence/extraction of the finite Zhu module and constituents is an interface.\\
Lemma 2.9 & \path{lemma_2_9_detection_from_graded_kernel_literature} & Lowest nonzero kernel degree and negative-wall arithmetic internal; irreducible DS nonvanishing and transfer are literature inputs.\\
Proposition 2.10 & \path{proposition_2_10_from_graded_kernel_literature} & Kernel contradiction, injectivity, and use of surjectivity internal from explicit DS exactness/detection hypotheses.\\
Theorem 2.11 & \path{theorem_2_11_from_graded_kernel_literature}; \path{minimal_reduction_maximality_core} & Injectivity-to-quotient-isomorphism and candidate-ideal equality internal.\\
Theorem 2.13 & \path{ds_ramond_casimir_maximality_core} & Common reduced-simplicity plus affine-lifting architecture kernel-checked in abstract form.\\
\midrule
Lemma 3.3--Prop. 3.6 & \path{DSeries.proposition_3_6_exact_images_from_literature_and_PBW}; \path{DSeries.proposition_3_6_exact_images_over_field} & Wall values and scalar-line deduction internal; DS nonvanishing/cyclicity and concrete PBW-space membership external.\\
Proposition 3.8 & \path{DSeries.proposition_3_8_exact_reduction} & Algebraic quotient isomorphism internal after exact-image and exact-sequence data are supplied.\\
Proposition 3.12 & \path{DSeries.weightedCasimir}; \path{DSeries.weightedCasimir_le}; trace/gap package & Finite Casimir arithmetic and downstream inequality checked; origin of the traced Ramond identity is boundary data.\\
Props. 3.17--3.20 & D-series spectral-bottom and ground-space packages & Lattice/inequality consequences are checked where encoded; concrete Li-twisted lowest-space realization is not rebuilt.\\
Theorem 3.21 & \path{DSeries.theorem_3_21_simplicity_core}; \path{DSeries.theorem_3_21_reduced_simplicity_certificate} & Casimir-gap contradiction and vacuum forcing internal from explicit Ramond interfaces.\\
Theorem 3.23 & \path{DSeries.theorem_3_23_gap_to_affine_graded_literature} & Reduced simplicity to affine maximality checked with graded DS exactness/nonvanishing exposed as inputs.\\
Proposition 3.27 & \path{DSeries.proposition_3_27_exact_images_from_literature_and_PBW}; \path{DSeries.proposition_3_27_exact_images_over_field} & Endpoint scalar-line conclusion internal from explicit PBW/DS inputs.\\
Proposition 3.28 & \path{DSeries.proposition_3_28_exact_reduction} & Endpoint quotient algebra internal once exactness/range data are supplied.\\
Propositions 3.31--3.32 & \path{DSeries.endpointA}; \path{endpoint_first_gap_formula}; endpoint ground package & Endpoint trace and first-shift gap arithmetic internal; collapsed-current/lowest-space interpretation external.\\
Theorem 3.33 & \path{DSeries.theorem_3_33_simplicity_core}; \path{DSeries.theorem_3_33_reduced_simplicity_certificate} & Endpoint reduced-simplicity contradiction checked.\\
Theorem 3.34 & \path{DSeries.theorem_3_34_gap_to_affine_graded_literature} & Endpoint affine kernel vanishing and maximal-ideal equality internal from graded DS inputs.\\
\midrule
Lemmas 4.1--4.2 / Prop. 4.3 & D4 extremal-PBW package; \path{D4.proposition_4_3_exact_images_over_field} & Triality wall/profile and extremal-square logic internal; concrete BRST nonvanishing and PBW realization external.\\
Proposition 4.4 & common exact-reduction quotient package & Exact quotient logic internal; DS exactness/cyclicity and actual reduced ideal interpretation supplied.\\
Lemmas 4.5--4.7 & \path{D4.AllowedTriple}; D4 Casimir/contraction package & Finite $A_1^3$ enumeration and Casimir arithmetic internal after the allowed-type and universal relation interfaces.\\
Proposition 4.8 & \path{D4.proposition_4_8_from_traced_relation}; \path{D4.proposition_4_8_matches_A} & Trace specialization and simplification checked by Lean.\\
Props. 4.10--4.12 & D4 PBW spectral-lattice / bottom degree-weight package & Spectral inequalities and finite eliminations internal in the model; actual Li-twist/Ramond interpretation remains a bridge.\\
Theorem 4.13 & \path{D4.d4_reduced_simplicity_certificate} & Triality Casimir-gap contradiction and vacuum conclusion kernel-checked.\\
Theorem 4.14 & \path{D4.d4_theorem_4_14_graded_literature} & Final affine lifting and candidate-ideal equality checked from explicit graded DS hypotheses.\\
\midrule
Theorem 5.6 / Prop. 5.7 & \path{E6.proposition_5_7_exact_images_from_literature_and_PBW}; \path{E6.proposition_5_7_exact_images_over_field} & Level/exponent and scalar-line deduction internal; singularity/DS/PBW realization inputs explicit.\\
Proposition 5.12 & \path{E6.proposition_5_12_exact_reduction} & Algebraic exact-reduction quotient checked from exactness, range, and exact-image data.\\
Proposition 5.14 & \path{E6.traceA}; \path{E6.trace_specialization_one}; \path{E6.trace_specialization_two} & Rational trace specialization checked; universal Ramond relation/contraction interpretation supplied.\\
Lemma 5.15 & \path{E6.Cmax}; \path{E6.Cnext}; pointwise Casimir bounds & Finite numerical maxima/second maxima used by the gap proof are checked after constituent restrictions are encoded.\\
Proposition 5.17 & E6 ground/spectral interface plus \path{positive_shift_gap} & Downstream integral-shift gap internal; actual Li-twisted spectrum and extremal line are representation-theoretic inputs.\\
Theorem 5.18 & \path{E6.theorem_5_18_simplicity_core}; \path{E6.theorem_5_18_reduced_simplicity_pointwise} & Pointwise Casimir bounds are converted to the averaged contradiction inside Lean.\\
Theorem 5.19 & \path{E6.theorem_5_19_pointwise_graded_literature} & Reduced simplicity and graded DS detection are assembled into affine maximality.\\
\midrule
Proposition 6.4 & \path{E7.section_6_2_exact_images_from_literature_and_PBW}; \path{E7.proposition_6_4_exact_images_over_field} & Exact current-power image logic checked for all three new $E_7$ cases from explicit DS/PBW inputs.\\
Proposition 6.5 & \path{E7.proposition_6_5_exact_reduction} & Algebraic exact-reduction quotient internal.\\
Proposition 6.10 / Cor. 6.11 & \path{E7.traceA}; \path{E7.trace_specializations} & Uniform trace formula specializations and rational arithmetic checked.\\
Props. 6.12--6.17 & E7 spectral/bottom package; \path{E7.Cmax}; \path{E7.Crest}; gap theorems & Low-degree exclusions and casewise Casimir inequalities encoded/checked; semantic constituent extraction remains boundary data.\\
Theorem 6.18 & \path{E7.theorem_6_18_simplicity_core}; \path{E7.theorem_6_18_reduced_simplicity_pointwise} & Positive-shift and ground-extremal contradictions assembled in Lean.\\
Theorem 6.19 & \path{E7.theorem_6_19_pointwise_graded_literature} & Final affine maximality checked from the reduced certificate and graded DS interfaces.\\
\midrule
Proposition 7.4 / Cor. 7.6 & \path{E8.corollary_7_6_exact_images_from_literature_and_PBW}; \path{E8.corollary_7_6_exact_images_over_field} & Wall arithmetic and scalar current-power conclusion internal; actual DS class/cyclicity/PBW realization supplied.\\
Proposition 7.7 & \path{E8.proposition_7_7_exact_reduction} & Algebraic exact-reduction quotient internal.\\
Proposition 7.14 & \path{E8.traceA}; \path{E8.ground_trace_formula}; \path{E8.first_gap_formula} & Uniform trace arithmetic and exact first-gap simplification checked.\\
Theorem 7.18 / Cor. 7.19 & \path{E8.Cmax}; \path{E8.Csecond}; \path{ground_minus_second_nonlevelone}; \path{ground_extremal_gap} & Finite Casimir maximum/second-gap consequences checked after allowed-weight data are encoded.\\
Proposition 7.22 & E8 spectral/ground interface plus \path{positive_shift_gap} & Numerical positive-shift consequence internal; concrete Li-twist spectrum is external.\\
Theorem 7.23 & \path{E8.theorem_7_23_simplicity_core}; \path{E8.theorem_7_23_reduced_simplicity_pointwise} & All five reduced-simplicity numerical branches assembled in Lean.\\
Theorem 7.24 & \path{E8.theorem_7_24_pointwise_graded_literature} & Final affine maximality checked from reduced simplicity and graded DS inputs.\\
\midrule
Proposition 8.1 & \path{RankReduction.proposition_8_1_exact_images_over_field} & Nonzero-on-a-line scalar extraction internal; concrete three BRST-image lines supplied.\\
Lemma 8.2 & \path{RankReduction.lemma_8_2_generation_core} & Abstract generation consequence checked; actual surviving-current generation fact supplied.\\
Lemma 8.3 & \path{RankReduction.lemma_8_3_quadratic_relation_core} & Algebraic elimination from the supplied quadratic relation checked.\\
Proposition 8.4 & \path{RankReduction.proposition_8_4_factor_core} & Factorization/surjection logic internal; concrete current-quotient realization external.\\
Theorem 8.5 & \path{RankReduction.theorem_8_5_induction_step_core}; \path{RankReduction.theorem_8_5_induction_step_graded_literature} & Sandwich/bijectivity, induction bookkeeping, and affine lifting internal; base case and rank-reduction bridge supplied.\\
\bottomrule
\end{longtable}
}

\section{Kernel-check evidence and axiom audit}

The verification claim rests on three distinct pieces of evidence, which should not be conflated.

\subsection{Successful project build}

The source frozen in GitHub release \code{v1.0.0} was built from the project root with
\begin{verbatim}
lake build
\end{verbatim}
and the supplied terminal transcript ends with
\begin{verbatim}
Build completed successfully (8708 jobs).
\end{verbatim}
Thus every declaration in the configured target elaborated and the generated proof terms were accepted in the pinned Lean/Mathlib environment.  The job count is not part of the mathematical claim and may change with cache state or Lake internals.

\subsection{Static source integrity and frozen snapshot}

A direct scan of the frozen file gives 14,296 lines and 446,445 bytes.  It contains 409 declarations introduced with the keyword \code{theorem}, 148 definitions introduced with \code{def}, nine \code{structure} declarations, and fifteen \code{inductive} declarations.  No declaration is introduced with \code{axiom}.  No executable proof contains a \code{sorry} or \code{admit} placeholder; the only textual occurrence of those words is in the explanatory boundary comment stating that such placeholders are absent.  The SHA-256 digest is
\begin{center}
\code{c4d73a23eeb58976bae7bd264b48a19a102e930576bb49a01b0c32d09f895906}.
\end{center}
This digest is also recorded in \code{SHA256SUMS.txt}, so the report can be tied to an exact source snapshot.

\paragraph{Frozen pre-build header note.}
The first comment block of the Lean source still contains the sentence that the statically merged source ``has not yet been certified here by a successful lake build.''  That sentence records the status at the time the source text itself was frozen.  In accordance with the release policy used here, the Lean file has not been edited merely to update that comment.  The later external verification run reported \code{Build completed successfully (8708 jobs)}, and that run is the build evidence asserted by this report.  Thus the stale header sentence should be read as historical metadata inside the immutable source snapshot, not as the current release status.

\subsection{Live \texttt{\#print axioms} targets}

The final lines of the Lean source execute sixteen live dependency queries.  They are grouped as follows.

{\footnotesize
\begin{longtable}{>{\raggedright\arraybackslash}p{0.22\textwidth} >{\raggedright\arraybackslash}p{0.68\textwidth}}
\toprule
Group & Targets\\
\midrule
\endfirsthead
\toprule
Group & Targets\\
\midrule
\endhead
Common cores & \path{ArakawaMoreau.minimal_reduction_maximality_core}; \path{ArakawaMoreau.ds_ramond_casimir_maximality_core}.\\
Affine maximality wrappers & \path{DSeries.theorem_3_23_gap_to_affine_graded_literature}; \path{DSeries.theorem_3_34_gap_to_affine_graded_literature}; \path{D4.d4_theorem_4_14_graded_literature}; \path{E6.theorem_5_19_pointwise_graded_literature}; \path{E7.theorem_6_19_pointwise_graded_literature}; \path{E8.theorem_7_24_pointwise_graded_literature}; \path{RankReduction.theorem_8_5_induction_step_graded_literature}.\\
Exact-image wrappers & \path{DSeries.proposition_3_6_exact_images_over_field}; \path{DSeries.proposition_3_27_exact_images_over_field}; \path{D4.proposition_4_3_exact_images_over_field}; \path{E6.proposition_5_7_exact_images_over_field}; \path{E7.proposition_6_4_exact_images_over_field}; \path{E8.corollary_7_6_exact_images_over_field}; \path{RankReduction.proposition_8_1_exact_images_over_field}.\\
\bottomrule
\end{longtable}
}

The terminal tail retained for this release visibly records the following representative outputs:
{\footnotesize
\begin{longtable}{>{\raggedright\arraybackslash}p{0.57\textwidth} >{\raggedright\arraybackslash}p{0.33\textwidth}}
\toprule
Target & Reported dependencies\\
\midrule
\endfirsthead
\toprule
Target & Reported dependencies\\
\midrule
\endhead
\path{D4.proposition_4_3_exact_images_over_field} & \code{propext}, \code{Classical.choice}, \code{Quot.sound}\\
\path{E6.proposition_5_7_exact_images_over_field} & \code{propext}, \code{Classical.choice}, \code{Quot.sound}\\
\path{E7.proposition_6_4_exact_images_over_field} & \code{propext}, \code{Classical.choice}, \code{Quot.sound}\\
\path{E8.corollary_7_6_exact_images_over_field} & \code{propext}, \code{Classical.choice}, \code{Quot.sound}\\
\path{RankReduction.proposition_8_1_exact_images_over_field} & \code{propext}, \code{Quot.sound}\\
\bottomrule
\end{longtable}
}
These are standard foundational Lean/Mathlib dependencies, not user-declared representation-theoretic axioms.  The report does not infer unrecorded exact output lists for the other live targets merely from build success; the source itself contains the commands so that a reproducing user can inspect the complete output directly.

\subsection{Why theorem parameters are different from axioms}

The absence of project-local axioms is not the same as a foundational formalization of the external mathematics.  For example, the exact-image wrappers quantify over hypotheses such as DS nonvanishing, exact transfer, cyclicity, PBW spanning, and concrete class membership.  Such hypotheses are ordinary theorem parameters.  Lean kernel-checks the implication from those hypotheses to the conclusion, while \code{\#print axioms} reports only logical constants appearing in the proof term.  Accordingly, the formal claim is conditional exactly where the theorem signatures say it is conditional.

\section{Reproducibility}

The release package pins the following environment:
\begin{itemize}
  \item Lean toolchain: \code{leanprover/lean4:v4.34.0-rc1};
  \item Mathlib tag requested by \code{lakefile.toml}: \code{v4.34.0-rc1};
  \item Mathlib commit resolved by \code{lake-manifest.json}:\\
  \code{de5ce8a9a66a4aa68a9bdbb35b63a06d34d9ca11}.
\end{itemize}

The recommended arXiv ancillary layout is
\begin{verbatim}
anc/
  ArakawaMoreauFormalAudit.lean
  README.md
  SHA256SUMS.txt
  lakefile.toml
  lake-manifest.json
  lean-toolchain
\end{verbatim}
The Lean file in this directory is the same frozen source whose SHA-256 digest is reported above; the report-writing process does not modify that file.

From the ancillary project directory, the primary verification command is
\begin{verbatim}
lake build
\end{verbatim}
The frozen \code{v1.0.0} release completed with
\begin{verbatim}
Build completed successfully (8708 jobs).
\end{verbatim}
The job count is cache- and implementation-dependent; successful termination without an error is the relevant reproducibility criterion.  A direct single-source check may additionally be run with
\begin{verbatim}
lake env lean ArakawaMoreauFormalAudit.lean
\end{verbatim}

For integrity checking, compare the source against \code{SHA256SUMS.txt}.  The expected digest for \code{ArakawaMoreauFormalAudit.lean} is
\begin{center}
\code{c4d73a23eeb58976bae7bd264b48a19a102e930576bb49a01b0c32d09f895906}.
\end{center}

\subsection{Release-consistency checklist}
The ancillary metadata have been normalized to the project name \code{ArakawaMoreauFormalAudit}: the main source is \code{ArakawaMoreauFormalAudit.lean}, the Lake target has the same name, and the project name in the distributed manifest has been updated accordingly.  The dependency revisions themselves are unchanged.  The Lean source is byte-for-byte identical to the successfully built frozen source, as witnessed by the digest above.

\section{What this verification does and does not claim}\label{sec:claim-boundary}

The strongest accurate summary of the project is:
\begin{quote}
The Lean bundle provides a kernel-checked formal audit of the paper-local proof architecture, finite arithmetic and Casimir computations, extremal PBW deductions, reduced-simplicity logic, and affine maximality lifting, conditional on explicitly stated VOA/BRST/DS and case-realization interfaces.
\end{quote}

In particular, the report makes three claims and only these claims: (i) the frozen source contains no project-local \code{axiom} declaration and no \code{sorry}/\code{admit} proof placeholder; (ii) the configured project has completed \code{lake build} successfully in the pinned environment; and (iii) the listed paper-local deductions are therefore checked by the Lean kernel from the hypotheses visible in their theorem signatures.  It does \emph{not} claim that the external VOA/BRST/DS hypotheses themselves have been derived from foundational definitions in Mathlib, nor does it claim an independent Lean reconstruction of the previously published exceptional $n=0$ cases used by the paper when assembling Theorem~1.1.

Accordingly, phrases such as ``fully formalized from definitions'' or ``unconditional Lean proof of all background representation theory'' would overstate the content of the bundle.  The intended description is ``kernel-checked formal audit modulo explicit literature/foundational interfaces.''

\paragraph{Result of the final report audit.}
The paper-number/declaration cross-check found no mismatch requiring a change to the frozen Lean source.  The principal correction made at the report level is scope clarification: the formal bundle audits the new case-level deductions and the alternative $D_\ell$ endpoint proof, while the already published exceptional $n=0$ cases remain literature inputs to the paper-level completion statement.  The report also records explicitly why the pre-build sentence in the immutable Lean header does not contradict the later successful release build.

\section{Availability of the formalization}

The formalization is distributed in two complementary forms.  The arXiv ancillary bundle accompanying this report contains the frozen Lean source, the Lake project metadata, the pinned toolchain, a concise README, and SHA-256 checksums.  A public maintainable copy is hosted at
\begin{center}
\url{https://github.com/jshemail12345-debug/ArakawaMoreauFormalAudit}.
\end{center}
The exact Lean source and Lake-metadata snapshot corresponding to this report is frozen as GitHub release \code{v1.0.0}:
\begin{center}
\url{https://github.com/jshemail12345-debug/ArakawaMoreauFormalAudit/releases/tag/v1.0.0}.
\end{center}
The arXiv ancillary copy and the tagged GitHub release therefore serve complementary reproducibility roles: the former is attached to the report version and may contain report-facing documentation, while the latter provides a fixed clonable public code snapshot.

\paragraph{Ancillary files.}
\begin{itemize}[nosep]
  \item \code{ArakawaMoreauFormalAudit.lean};
  \item \code{README.md} and \code{SHA256SUMS.txt};
  \item \code{lean-toolchain}, \code{lakefile.toml}, and \code{lake-manifest.json}.
\end{itemize}

The GitHub release also provides the bundled archive \code{ArakawaMoreauFormalAudit\_release.zip} in addition to GitHub's automatically generated source archives.

\section*{Acknowledgements}
The author acknowledges the assistance of the language model DeepSeek-V4-Pro in drafting and refining portions of the Lean 4 formalization and in preparing this verification report.  The author remains responsible for the mathematical statements, the formalization boundary, the verification claims, and the final contents of the released source and report.

\end{document}